\documentclass[letterpaper, 10 pt, conference]{ieeeconf}  % Comment this line out
\IEEEoverridecommandlockouts                              % This command is only
\usepackage{graphics} % for pdf, bitmapped graphics files
\usepackage{epsfig} % for postscript graphics files
\usepackage{mathptmx} % assumes new font selection scheme installed
\usepackage{amsmath} % assumes amsmath package installed
\usepackage{amssymb}  % assumes amsmath package installed
\usepackage{color}
\usepackage{float}
\usepackage{algorithmic}
\usepackage{algorithm}[0.1]
\usepackage{subfigure}
\usepackage{caption}

\newcommand{\ben}{\begin{equation}}
\newcommand{\een}{\end{equation}}

\newtheorem{thm}{Theorem}[section]

\newtheorem{dfn}[thm]{Definition}
\newtheorem{eg}[thm]{Example}

\newcommand{\eps}{\varepsilon}

\newcommand{\RR}{\mathbb{R}}

\newcommand{\inn}[2]{\left\langle #1, #2 \right\rangle} 
\newcommand{\norm}[2]{\left\| #1 \right\|_{#2} }

\newcommand{\overbar}[1]{\mkern 1.5mu\overline{\mkern-1.5mu#1\mkern-1.5mu}\mkern 1.5mu}

\newcommand{\DD}{\mathsf{\Omega}}

\newcommand{\beqn}{\begin{eqnarray*}}
	\newcommand{\eeqn}{\end{eqnarray*}}
\newcommand{\beq}{\begin{eqnarray}}
\newcommand{\eeq}{\end{eqnarray}}
\newcommand{\benn}{\begin{equation*}}
\newcommand{\eenn}{\end{equation*}}

\title{\LARGE \bf
Optimal Actuator Design for Vibration Control Based on LQR Performance and Shape Calculus
}

\author{M. Sajjad Edalatzadeh, Dante Kalise, Kirsten A. Morris and Kevin Sturm% <-this % stops a space
\thanks{Sajjad Edalatzadeh is with Department of Applied Mathematics, University of Waterloo, Waterloo, ON N2L 3G1, Canada.
        {\tt\small email: msedalatzadeh@uwaterloo.ca}}%
\thanks{Dante Kalise is with Department of Mathematics, Imperial College London, London SW7 2AZ, UK.
	{\tt\small email: d.kalise-balza@imperial.ac.uk }}%}%
\thanks{Kirsten A. Morris is with Faculty of Mathematics, Applied Mathematics, University of Waterloo, Waterloo, ON N2L 3G1, Canada.
	{\tt\small email: kmorris@uwaterloo.ca}}%
\thanks{Kevin Sturm is with Institute for Analysis and Scientific Computing, Technical University of Vienna Wiedner Hauptstra\ss e 8-10, 1040 Wien, Austria
	{\tt\small email: kevin.sturm@tuwien.ac.at}}%
}

\begin{document}

\maketitle
\thispagestyle{empty}
\pagestyle{empty}

%%%%%%%%%%%%%%%%%%%%%%%%%%%%%%%%%%%%%%%%%%%%%%%%%%%%%%%%%%%%%%%%%%%%%%%%%%%%%%%%
\begin{abstract}
Optimal actuator design  for a vibration control problem is calculated. The actuator shape is optimized according to the closed-loop performance of the resulting linear-quadratic regulator and a penalty on the actuator size. The optimal actuator shape  is found by means of shape calculus and a topological derivative of the  linear-quadratic regulator (LQR) performance index. An abstract framework is proposed based on the theory for infinite-dimensional optimization of both the actuator shape and the associated control problem. A numerical realization of the optimality condition is presented for the actuator shape using a level-set method for topological derivatives. A Numerical example illustrating the design of actuator for Euler-Bernoulli beam model is provided.
\end{abstract}

%%%%%%%%%%%%%%%%%%%%%%%%%%%%%%%%%%%%%%%%%%%%%%%%%%%%%%%%%%%%%%%%%%%%%%%%%%%%%%%%
\section{INTRODUCTION}
Actuator location and shape are important design variables for  feedback synthesis  in active vibration control. Optimal actuator design improves performance of the controller and significantly reduces the cost of control. Optimized static controls can yield better performance on a structure than dynamic controls with actuation at different locations \cite{FahDem,morris2015comparison,Mnoise}. In the engineering literature, the optimal actuator location problem has gained considerable interest, as shown in \cite{frecker2003recent,van2001review,kubrusly1985sensors,valadkhan2010stability}. 

The optimal actuator location theory of linear distributed parameter systems has been extensively studied in the literature. In \cite{morris2011linear}, it was proven that an LQ-optimal actuator location exists if the  control operator continuously depends on actuator locations. Similar results have been obtained for $H_2$ and $H_\infty$ controller design objectives \cite{DM2013,kasinathan2013h}. Other objectives such as maximizing controllability  or stability margin have been studied in \cite{demetriou2000numerical,HebrardHenrot2005}.

In contrast, optimal actuator shape design has only been studied in relatively few works.  In \cite{PTZ2013}, the optimal shape and position of an actuator for the wave equation in one spatial dimension are discussed. An actuator is placed on a subset $\omega\in [0,\pi]$ with a constant Lebesgue measure $L\pi$ for some $L\in (0,1)$. The optimal actuator minimizes the norm of a Hilbert Uniqueness Method (HUM)-based control; such control steers the system from a given initial state to zero state in finite time. In \cite{privat2017actuator}, optimal actuator shape and position for  linear parabolic systems are discussed. This paper adopts the same approach as in \cite{PTZ2013} but with initial conditions that have randomized Fourier coefficients. The cost is defined as the average of the norm of HUM-based controls. In \cite{kalise2017optimal}, optimal actuator design for linear diffusion equations has been discussed. A quadratic cost function is considered, and shape and topological derivative of this function are derived. Numerical results show significant improvement in the cost and performance of the control. Optimal sensor design problem is in many ways similar to the optimal actuator design problem. In \cite{privat2015optimal}, optimal sensor shape design has been studied where the observability is being maximized over all admissible sensor shapes. Optimal actuator design problems for nonlinear distributed parameter systems has also been studied. In \cite{edalatzadehSICON}, it is shown that under certain conditions on the nonlinearity and the cost function, an   optimal input and actuator design exist, and optimality equations are derived. Results  are applied to the nonlinear railway track model as well as to the semi-linear wave model in two spatial dimensions. Numerical techniques to calculate the optimal actuator shape design are mostly limited to linear quadratic regulator problems, see eg.  \cite{allaire2010long,kumar1978optimal-a,kubrusly1985sensors,darivandi2013algorithm}. For controllability-based approaches,  numerical schemes have been studied in \cite{munch2011optimal,munch2009optimal,munch2013numerical}.

In this paper,  optimal actuator design  for the Euler-Bernoulli beam model arising in active vibration control is studied. We evaluate the performance of a given actuator in terms of the cost associated to the resulting LQR problem and a penalty on the actuator size. This defines a cost function for the actuator and corresponding optimality conditions. The optimal actuator design is found based on shape calculus and the notion of topological derivative. Well-posedness of the abstract problem for a linear evolution equation, and the particular setting for the Euler-Bernoulli model is established. The characterization of the optimal actuator shape by means of the topological derivative dictates the construction of an approximation scheme based on the discretization of the LQR problem, along with a level-set method for the realization of the optimality condition for the actuator shape.

The rest of the paper is structured as follows. In Section II, an abstract formulation of the optimal actuator design problem is described. Section III introduces the control setting associated to the the Euler-Bernoulli model for beam vibration control. Section IV, presents a numerical method for the approximation of optimal actuators based on level-set method for topological derivatives. In Section V, computational results are discussed.

\section{Abstract Formulation of the Optimal Actuator Design Problem}

Consider an abstract evolution equation over $[0,\tau]$
\begin{equation}\label{con-sys}
\begin{cases}
\dot{z}(t)=Az(t)+B(r)u(t),\quad t\in (0,\tau],\\
z(0)=z_0.
\end{cases}
\end{equation}
Let $H$ and $U$ be Hilbert spaces, $A:D(A)\to H$ be the state operator, and $z(t)$  the state of the system. The state operator $A$ generates a strongly continuous semi-group on $H.$ 
Let $r$ indicate an actuator design parameter taking values in a topological space $K$. The input operator $B(r):U\to H$ continuously depends on actuator design $r$, that is
\begin{equation}\label{cont}
\lim_{r_2\to r_1}\norm{B(r_2)-B(r_2)}{\mathcal L(U,H)}=0.
\end{equation}
The input $u$ is assumed to belong to the space $L^2(0,\tau;U)$.

Depending on the regularity of inputs $u$ and initial conditions $z_0$, various notions of solutions can be defined for this system  \cite[Definition II.3.2]{bensoussan2015book}. \begin{dfn}(Weak solution)
Let $A^*:D(A^*)\to H$ be the dual operator of $A$. If $z$ belongs to $L^2(0,\tau;H)$, and for all $h \in D(A^*)$, $\inn{h}{z(t)}$ belongs to $H^1(0,\tau)$, and for almost all $t \in [0, \tau]$ and $h$ in $D(A^*)$
\begin{equation}\label{weak}
\begin{cases}
\frac{d}{dt}\inn{h}{z(t)}=\inn{A^*h}{z(t)}+\inn{h}{B(r)u(t)},\\
z(0)=z_0,
\end{cases}
\end{equation}
then $z(t)$ is said to be a {\it weak solution} of (\ref{con-sys}).
\end{dfn}
Since the  state space $H$ is a Hilbert space, this definition of weak solution is equivalent to a weak solution defined in the variational framework for parabolic systems \cite[Remark II.3.2]{bensoussan2015book}.

Let $U_{ad}\subset L^2(0,\tau;U)$ be the set of admissible inputs. The admissible actuator shape set is denoted by $K_{ad}$ and is compact with respect to the topology on $K$. Given an admissible actuator parameter $r$, define the  finite-horizon  cost 
\begin{align} 
J (u,r;z_0):=\frac{1}{2}\int_0^\tau \|z(t)\|^2_{H}+\gamma\|u(t)\|^2_{U}\,dt.
\end{align}
for some positive number $\gamma$, constrained to the evolution \eqref{con-sys}. The actuator shape will be designed based on minimizing $J (u,r;z_0)$ over all admissible inputs and actuator shapes, that is,
\begin{align}
 \underset{r\in K_{ad}}{\min}\;\underset{u\in U_{ad}}{\min}&J  (u,r;z_0)\,.
\end{align}
 
Provided that $K_{ad}$ is compact in the topology on $K,$ Theorem 4.1 in \cite{edalatzadehSICON} guarantees the existence of an optimal control and actuator shape for this problem. 

Let $\DD$ be a bounded set in $\mathbb{R}^n$ related to the spatial domain of the abstract evolution \eqref{con-sys}. In order to characterize the optimal shape in terms of 
\begin{align}
J_1(r,z_0):=\underset{u\in U_{ad}}{\min}&J  (u,r;z_0)\,,
\end{align}
Analysis is restricted  to the case in which $r$ corresponds to a subset $\omega$ of $\Omega$. In this setting, we define the \textsl{topological derivative} of $J_1(\omega)$, where we omit the dependence with respect to $z_0$.
\begin{dfn}[Topological derivative]
 	The topological derivative of $J_1$ at $\omega$ in the point $\eta_0\in \DD \setminus \partial \omega$ is defined by
 	\begin{equation}\notag
 	\mathcal T J_1(\omega)(\eta_0) =\begin{cases}
 		\lim_{\eps\searrow 0}\frac{J_1(\omega\setminus \bar B_\eps(\eta_0)) - J_1(\omega)}{|\bar B_\eps(\eta_0)|} & \text{ if } \eta_0 \in \omega, \\[2mm]
 		\lim_{\eps\searrow 0}\frac{J_1(\omega\cup  B_\eps(\eta_0)) - J_1(\omega)}{|B_\eps(\eta_0)|} & \text{ if } \eta_0 \in \DD \setminus \overline \omega.
\end{cases}
 	\end{equation}
\end{dfn}

\section{Beam Vibrations}
Consider a simply supported Euler-Bernoulli beam with Kelvin-Voigt damping coefficient $C_d$ and viscous damping coefficient $\mu$.  The state $w(x,t)$ corresponds to the transverse vibrations of the beam  with $(x,t)\in [0,1]\times[0,\tau]$. The actuator shape is characterized by an indicator function $\chi_{\omega}(x)$, where  $\omega$ is a Lebesgue measurable subset of $[0,1]$.  The beam vibrations are governed by
\begin{equation}\notag
\begin{cases}
\begin{aligned}
&\partial^2_{tt} w+\partial^2_{xx}\left(\partial^2_{xx}w+C_d\partial^3_{xxt}w\right)+\mu \partial_tw=\chi_{\omega}(x)u(t),  &t&>0,\\
&w(x,0)=w_0(x),\quad  \partial_tw(x,0)=v_0(x),\\
&w(0,t)=w(1,t)=0,&t&\ge 0,\\
&\partial^2_{xx}w(0,t)+C_d\partial^3_{xxt}w(0,t)=0,   &t&\ge 0,\\
&\partial^2_{xx}w(1,t)+C_d\partial^3_{xxt}w(1,t)=0,  &t&\ge 0,
\end{aligned}
\end{cases}
\end{equation}
where the subscript $\partial_x$ denotes the derivative with respect to $x$; the derivative with respect to $t$ is indicated similarly. The velocity of the displacement is denoted by $v=\partial_t w$. We define the augmented state $z:=(w,v)$. Also, define the state space $H:=H^2(0,1 )\cap H_0^1(0,1 )\times L^2(0,1 )$  with norm
\begin{equation}
\| (w,v) \|_H^2=\int_0^{1} \left(\partial^2_{xx}w(x)\right)^2+ w^2(x)+v^2(x) \, dx \label{eq: norm},
\end{equation}
and  the closed self-adjoint positive operator 
\begin{equation}\label{eqn-Ao}
\begin{cases}
A_0w:=\partial^4_{xxxx}w, \\
D(A_0):=\left\lbrace w\in H^4(0,1 )| \, w(0)=w(1)=0,\right.\\
\qquad \qquad \; \left. \partial^2_{xx}w(0)=\partial^2_{xx}w(1)=0 \right\rbrace.
\end{cases}
\end{equation}
The state operator is 
\begin{flalign*}
A(w,v) &=\left(v,-A_0(w+C_dv)-\mu v\right),\\
D(A)=&\lbrace(w,v)\in H| \, v\in H^2(0,1 )\cap H_0^1(0,1 ),   \\
&\quad \quad
w+C_dv\in D(A_0)  \rbrace.
\end{flalign*}
Letting $U=\RR$ and $U_{ad}=L^2(0,\tau),$   for $r= \chi_{\omega},$ define the input operator 
\begin{equation}
B(r)u=(0,\chi_{\omega}  u).
\end{equation}

Define $K:=L^\infty(0,1)$. The admissible actuator design set is for some $V>0$ and $0<c<1$, 
\begin{flalign*}
K_{ad}:=\lbrace\chi_\omega\in BV(0,1) :  Var\{\chi_\omega(x),[0,1]\}\le V, \; |\omega| = c \rbrace.
\end{flalign*}
where  $BV(0,1)$ denotes the set of functions of bounded variations on $(0,1)$, $Var\{r,[0,1]\}$ denotes the total variation of $r$ on $[0,1]$. Since $BV(0,1)$ is compactly embedded into $L^1(0,1)$ 
and since $r\mapsto Var(r,[0,1])$ is lower semi-continuous, it follows that $K_{ad}$ is compact in $K$ with respect to weak-* topology; see \cite[Thm. 1.126]{barbu2012convexity}. 
Also, for any sequence in $K_{ad}$, there is a subsequence that converges point-wise in $(0,1)$ to a function in $K_{ad}.$ Thus, if $r_2\in K_{ad}$ converges to $r_1\in K_{ad}$ in weak-* topology, so does it point-wise. Noting that
\begin{equation}
\norm{B(r_2)-B(r_1)}{\mathcal L(U,H)}=  \left( \int_0^1 (r_2 (x) - r_1(x) )^2 dx\right)^{\frac{1}{2}},
\end{equation}
condition \eqref{cont} follows from this and the Dominated Convergence Theorem.

For every initial condition $z_0=(w_0,v_0)$, define the cost
\begin{align}\notag
J_1(\omega, z_0) := \min_{u} \int_0^\tau \int_0^1 \norm{(w,v)}{}^2 + \gamma|u(t)|^2\,dt.
\end{align}
An optimal shape and control for this problem exists \cite[Thm. 4.1]{edalatzadehSICON}. 
For given $z_0$ and $\omega$, the associated optimal control will be denoted by $\bar u^{\omega,z_0}$.

The weak solution to the beam equation is obtained using variational calculus. That is, the weak solution satisfies 
\begin{align}
\int_0^\tau \int_0^1\partial^2_{tt} w\varphi &+ (\partial^2_{xx}w+C_d\partial^3_{xx t}w) \partial^2_{xx} \varphi + \mu\partial_t w\varphi \; dx dt \notag \\
&=\int_0^\tau\int_0^1 \chi_{\omega} u(t)\varphi \;dxdt. \label{sys}
\end{align}

 The corresponding optimal deflection associated with the optimal control $\bar u^{\omega,z_0}$ is denoted by $\bar w^{\omega,z_0}$. Then the topological derivative of $J_1$ at an open set $\omega$ in  $\eta_0 \in \DD\setminus \partial \omega$ is given by $\mathcal T J_1(\omega)(\eta_0) =$
\begin{align}
 		 \left.
 		\begin{cases}
 			-\int_0^T\overline u^{\omega,z_0}(\eta_0,s)\bar p^{\omega,z_0}(\eta_0,s)\;ds & \text{ if } \eta_0 \in \omega, \\
 			\int_0^T\overline u^{\omega,z_0}(\eta_0,s)\bar p^{\omega,z_0}(\eta_0,s)\;ds &  \text{ if }  \eta_0 \in  K_{ad}\setminus \overline \omega,
 		\end{cases}
 		\right.
\end{align}
where the adjoint $\bar p^{\omega,z_0}$ satisfies
\begin{align}
&\int_0^\tau \int_0^1\partial^2_{tt} \psi \bar p^{\omega,z_0} + (\partial^2_{xx}\psi+C_d\partial^3_{xx t}\psi) \partial^2_{xx}\bar p^{\omega,z_0} + \mu\partial_t \psi \bar p^{\omega,z_0} dxdt \notag\\
&\quad  = -\int_0^\tau\int_0^1 2 (\bar w^{\omega,z_0} + \partial_t\bar{w}^{\omega,z_0})\psi  \;dxdt.\label{adj}
\end{align}

\section{Numerical Approximation}

In this section, a finite-dimensional approximation of the state and adjoint variables is introduced, and the subsequent derivation of an approximation to the optimal actuator design problem. Assume $H$ is a separable Hilbert space. Let $V_N$ be an $N$-dimensional subspace of $D(A^*)$, and $\{\phi_n\}_{n=1}^N$ be an orthonormal  basis for $V_N$. An approximate state $z_N$ is defined as
\begin{equation}
z_N(x,t):=\sum_{n=1}^N\phi_n(x)a_n(t).
\end{equation}
At time $t=0$, the states $a_n(0)$ satisfy  
\begin{equation*}
a_n(0)=\inn{z_0}{\phi_n}_{H}.
\end{equation*}
Using test functions $h=\phi_n$ in \eqref{weak}, letting $Z_N=[a_1(t),a_2(t),\cdots,a_N(t)]$, and $U:=\mathbb{R}$ result in the finite-dimensional system
\begin{equation}
\dot{Z}_N(t)=A_NZ_N(t)+B_Nu(t).
\end{equation}
where $A_N\in \mathbb{R}^{N\times N}$ and $B_N \in \mathbb{R}^{N\times 1}$ are matrices with entries 
\begin{flalign*}
(A_N)_{nm}=\inn{A^*\phi_n}{\phi_m}, \quad (B_N)_{n}=\inn{\phi_n}{B(r)}.
\end{flalign*}
The following discrete quadratic cost function which includes a penalty for the actuator size $|\omega|$ is considered
\begin{flalign}\label{eq:min_problem_discrete}
J_N(u_N,\omega;Z_N(0)):=&\frac{1}{2}\int_0^\tau \|z_N(\cdot,t)\|^2_{H}+\gamma\|u(t)\|^2_{U}\,dt\notag \\
&+\alpha \left(|\omega|-c\right)^2\,,
\end{flalign} 
and in analogy to the continuous problem
\begin{equation}
J_{1,N}(\omega;Z_N(0)):=\min_{u \in U_{ad}}J_N(u,\omega;Z_N(0)).
\end{equation}
The input $\bar{u}$ minimizing $J_N$ is obtained by solving the  differential Riccati equation
\begin{equation}\notag
\begin{cases}
\dot{\Pi}_N(t)=-A^T_N\Pi_N(t)-\Pi_N(t)A_N+\frac{1}{\gamma}\Pi_N(t)B_NB_N^T\Pi_N(t)\\
\qquad \qquad-I_N, \; \forall t\in [0,\tau),\\
\Pi_N(\tau)=0.
\end{cases}
\end{equation}
After solving for $\Pi_N(t)$, the optimal cost is evaluated as 
$$J_{1,N}(\omega;Z_N(0))=Z_N^T(0)\Pi_N(0)Z_N(0)+\alpha \left(|\omega|-c\right)^2.$$
The optimal input $\bar{u}(t)$ is
\begin{equation}
\bar{u}(t)=-\frac{1}{\gamma}B^T_N\Pi_N(t)\bar{Z}_N(t)\,,
\end{equation}
and the adjoint state associated with $\bar{Z}_N$ and $\bar{u}$ is 
\begin{equation}
\bar{p}_N(x,t)=\sum_{n=1}^N\phi_n(x)b_n(t),
\end{equation}
where $b_n(t)$ are the entries of $2\Pi_N(t)\bar{Z}_N(t)$.

We obtain a numerical realization of the topological derivative associated to the optimal actuator design problem. Given an actuator $\omega$, define the function 
\begin{equation}
g_N^{\omega}(x) =
-\int_0^\tau\bar{u}(t)\bar{p}_N(x,t)\;dt + 2\alpha(|\omega|-c).
\end{equation}
The necessary optimality condition for an actuator shape $\omega$ is
\ben\label{eq:necessary_g}
\begin{split}
	g_N^{\omega}(x) & \le 0 \quad \text{ for all } x\in  \omega,\\
	g_N^{\omega}(x) & \ge 0 \quad \text{ for all } x\in \DD\setminus \omega.
\end{split}
\een

In order to find a stationary point for the topological derivative, we use a level-set method as proposed in  \cite{AmstutzHeiko06}. A remarkable advantage of this formulation is that the level-set evolution is versatile enough to allow the generation of multi-component actuators. We characterize the actuator $\omega$ by means of a level-set function  $\psi_N \in V_N$, such that 
\ben\label{eq:necessary_g}
\begin{split}
	\psi_N(x)  & <0 \quad \text{ for all } x\in  \omega,\\
	\psi_N(x)  & = 0 \quad \text{ for all } x\in \partial \omega.\\
	\psi_N(x)  & > 0 \quad \text{ for all } x\in \DD\setminus \bar\omega\,.\\
\end{split}
\een
Given an initial guess for $\psi_N^0$, the optimal actuator shape is determined by the iteration
\ben\notag
\psi^{i+1}_N = (1-\beta_i)\psi^i_N + \beta_i \frac{g_N^{\omega_i}}{\|g_N^{\omega_i}\|},\quad \omega_i := \{ x\in\DD: \psi^i_N(x)<0\},
\een
where $\beta_i$ is the step size of the level-set update, which evolves according to the realization of the topological derivative for the current shape. After reaching a stationary point of the topological derivative, the level-set evolution stops. Algorithm \ref{alg:topo} summarizes all the steps for the calculation of the optimal actuator shape. 
\begin{algorithm}[H]
	\begin{algorithmic}
		\STATE{\textbf{Input}: $\psi^0_N\in V_N(\DD)$, $\omega_0:=\{x\in \overbar{\DD},\psi^0_N(x)<0\}$, $\beta_0>0$, $Z_N(0)$, and tolerance $\eps >0$. }
		
		\WHILE{ $\|\omega_{i+1}-\omega_{i}\|\ge  \eps$ }
		\IF{ $ J_{1,N}(\omega_{i+1},Z_N(0)) < J_{1,N}(\omega_{i},Z_N(0))$}
		\STATE{ $\psi^{i+1}_N \gets (1-\beta_i) \psi^i_N +  \beta_i \frac{g_N^{\omega_i}}{\|g_N^{\omega_i}\|}$\;}
		\STATE{	$\beta_{i+1} \gets \beta_i$}
		\STATE{	$\omega_{i+1} \gets \{\psi^{i+1}_N<0\}$}
		\STATE{	$i\gets i+1$ }
		\ELSE
		\STATE{decrease $\beta_i$\;}
		\ENDIF
		\ENDWHILE
		\RETURN{ optimal actuator $\omega_{opt}$ }
		\caption{Level set algorithm for optimal actuator design}\label{alg:topo}
	\end{algorithmic}
\end{algorithm}
Algorithm \ref{alg:topo} also includes a line search strategy for finding a suitable $\beta_i$ to ensure a decay on the overall cost $J_{1,N}$. In practice, the algorithm is also embedded inside a continuation approach over the quadratic penalty parameter $\alpha$ in \eqref{eq:min_problem_discrete} to discard suboptimal actuators which satisfy the size constraint but which lead to poor closed-loop performance.

\section{Numerical Experiments}

\begin{figure}[!ht]
	\centering
	\subfigure[evolution of the  cost $J_{1,N}$  under the level-set set Algorithm \ref{alg:topo}. \label{J1a}]{\includegraphics[width=0.5\textwidth]{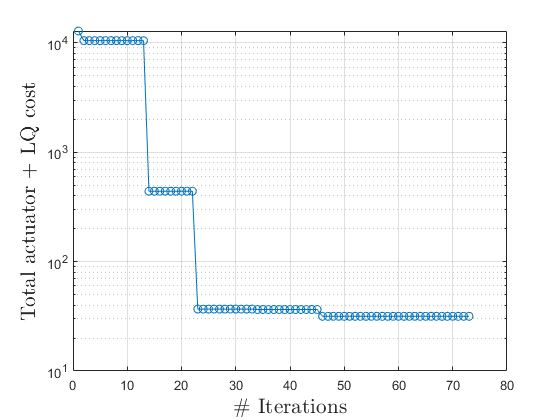}}
	\subfigure[optimal actuator \label{J1b}]{\includegraphics[width=0.5\textwidth]{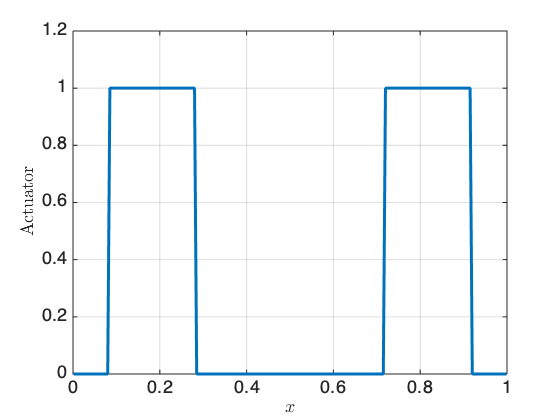}}
	\caption{\footnotesize Algorithm performance for the initial condition $w(x,0)=\sin(3\pi x)$ and $v(x,0)=0$ with  an actuator volume constraint of $40\%$ of the domain.}\label{J1}
\end{figure}

\begin{figure}[H]
	\centering
	\subfigure[uncontrolled displacement]{\includegraphics[width=0.4\textwidth]{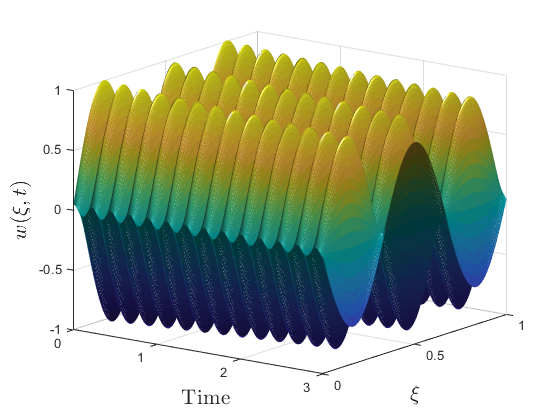}}
	\subfigure[controlled displacement]{\includegraphics[width=0.4\textwidth]{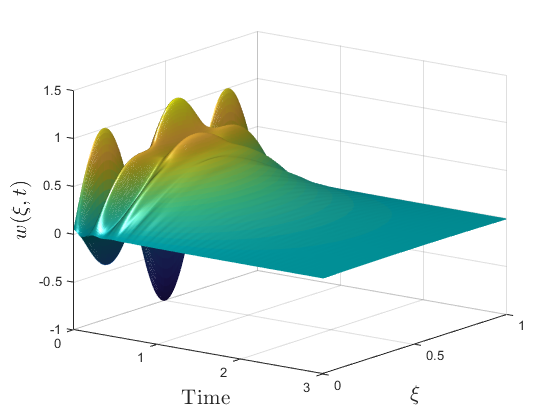}}
	\subfigure[controlled velocity]{\includegraphics[width=0.4\textwidth]{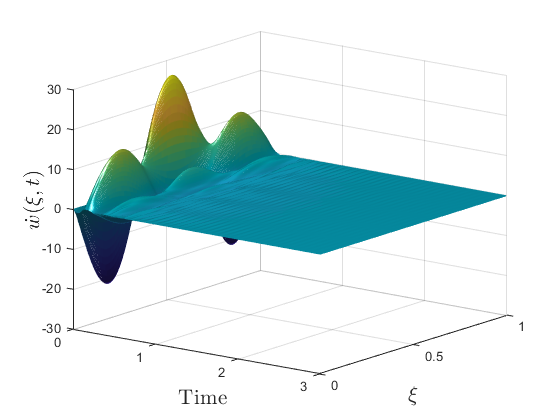}}
	\subfigure[optimal controller]{\includegraphics[width=0.4\textwidth]{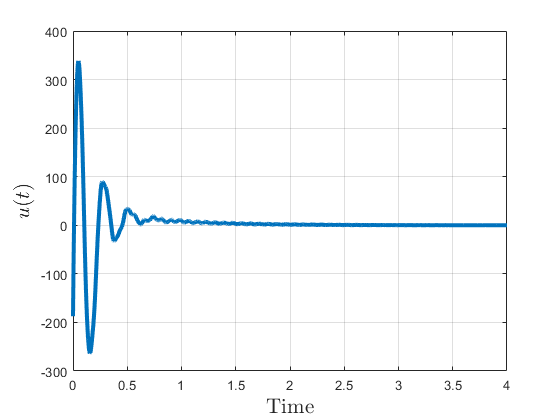}}
		\caption{\footnotesize Closed-loop performance of the optimal actuator.}\label{clopt}
\end{figure}

\begin{figure}[!ht]
	\centering
	\subfigure[controlled displacement]{\includegraphics[width=0.45\textwidth,height=0.2\textheight]{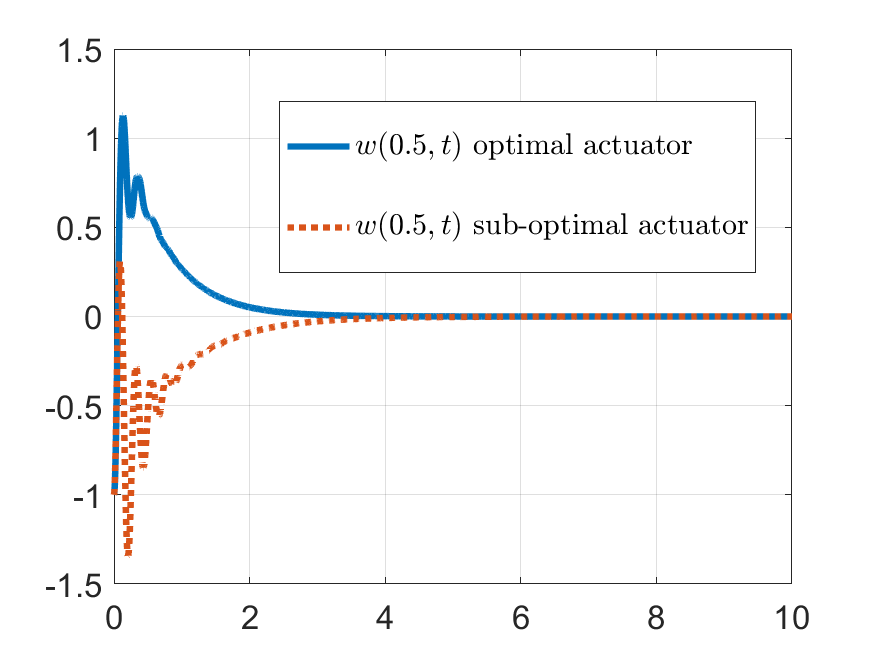}}
	\subfigure[controlled velocity]{\includegraphics[width=0.45\textwidth,height=0.2\textheight]{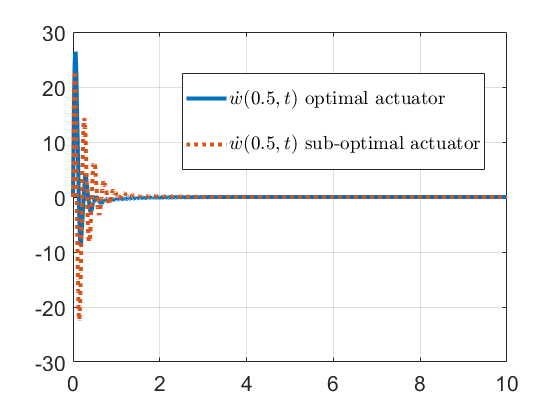}}
	\subfigure[optimal controller]{\includegraphics[width=0.45\textwidth,height=0.15\textheight]{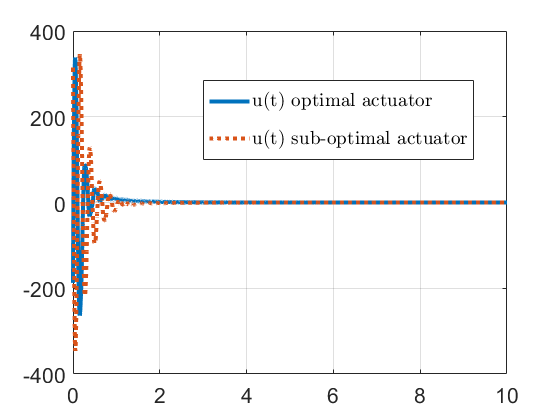}}
	\caption{\footnotesize Closed-loop performance of the optimal actuator against a sub-optimal actuator $w_s=[0.2,0.6]$ of the same volume at $x=0.5$. With the optimal actuator, better performance is achieved with a smaller control signal.  }\label{subopt}
\end{figure}

We present a series of numerical experiments illustrating the implementation and performance of the proposed method for actuator design. Let $x\in \DD:=(0,1)$, we consider the initial condition given by $w(x,0)=\sin(3\pi x)$ and $v(x,0)=0$, and the damping parameters set to $C_d=10^{-4}$ and $\mu=10^{-3}$. The system dynamics are discretized by taking the first $N=40$ eigenmodes of the operator $A$. The time horizon is set to $\tau=200$, the control penalty is set to $\gamma=10^{-3}$ and the volume constraint $c=0.4$, i.e. the actuator is enforced to cover $40\%$ of the domain. The continuation with respect to the volume penalty starts the execution of Algorithm \ref{alg:topo} with $\alpha=0.1$, increasing $\alpha$ by a factor of 10 until reaching $\alpha=10^4$. The level-set algorithm is stopped after a tolerance of $\epsilon=10^{-7}$ is reached. Every 20 iterations the level-set method is reinitialized with a signed distance function computed as in \cite{AFK15}.
The performance of the level-set algorithm is first shown in Figure \ref{J1a}, where we see the decrease of the overall cost $J_{1,N}$ as the shape of the actuator evolves. The initial actuator corresponds to $\omega_0=[0.1,0.9]$. The evolution is monotone as ensured by our line search step, and reaches the desired tolerance after 70 iterations, reducing the total cost in three orders of magnitude. The optimal actuator is shown in Figure \ref{J1b}, and we observe a splitting into two components, consistent with the shape of the initial condition under which the actuator was optimized. The closed-loop performance with the optimal actuator is shown in Figure \ref{clopt}, and it is compared against the closed-loop performance of a single-component, sub-optimal actuator $w_s=[0.2,0.6]$ in Figure \ref{subopt}.

\begin{figure}[!ht]
	\centering
	\subfigure[Viscous and Kelvin-Voigt damping]{\includegraphics[width=0.45\textwidth,height=0.2\textheight]{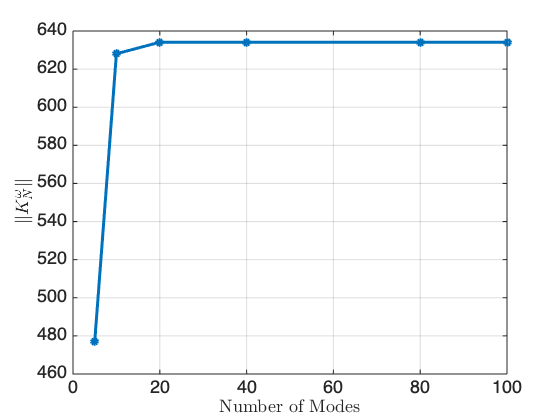}}
	\subfigure[Only viscous damping]{\includegraphics[width=0.45\textwidth,height=0.2\textheight]{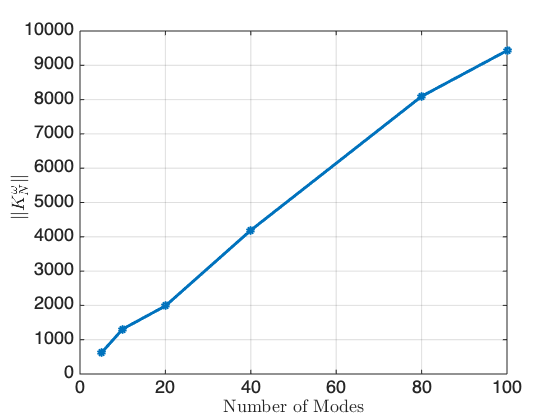}}
	\caption{\footnotesize Convergence of the Kalman gain and associated optimal actuator with resoect to the number of discretization modes depends on the choice of damping parameters.  }\label{kgood}
\end{figure}

\begin{figure}[!ht]
\centering
	\includegraphics[width=0.45\textwidth]{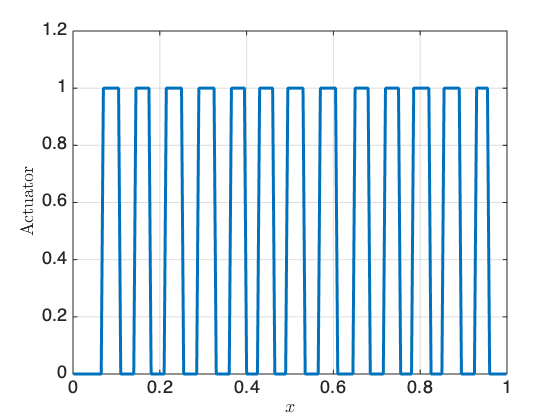}\hfill
	\caption{\footnotesize In the absence of Kelvin-Voigt damping, the optimal actuator fails to converge, and splits into multiple components as the number of discretization modes increases. In this figure, $N=40$ is the order of approximation. }\label{actbad}
\end{figure}

The convergence of the optimal shape with respect to the discretization of the system dynamics was investigated. In the presence of both viscous and Kelvin-Voigt damping, the optimal actuator design problem yields a numerical solution. In Figure \ref{kgood}(a), we illustrate this by showing the norm of the Kalman gain $K_N^{\omega}=-\frac{1}{\gamma}B^T_N\Pi_N(0)$ for increasing number of modes in our spectral discretization of the Euler-Bernoulli dynamics. In contrast, when  $C_d=0$, the semigroup generated by $A$ is no longer analytic, and the Kalman gain fails to converge with respect to the discretization, as shown in Figure \ref{kgood}(b). This translates into the lack of convergence for the actuator shape, which continues to split into more components as the number of modes increases, see Figure \ref{actbad}.

\subsection*{Conclusion and Future Research} 
We have discussed the optimal actuator design problem based on the LQR performance of the closed-loop associated to a given actuator shape. We presented a methodology to find the optimal actuator as a stationary point of the topological derivative of the cost associated to the LQ cost. We have applied this technique to solve the optimal actuator problem for a linear Euler-Bernoulli beam model.  The improvement in the performance of the system as well as convergence with respect to discretization parameters were illustrated. It was shown that without the presence of Kelvin-Voigt damping, the approximate optimal actuator shapes do not converge as the order of approximation is increased. Future work will focus on the optimal actuator shape design for nonlinear models of flexible structures. 

\bibliographystyle{ieeetr}
\bibliography{library}

\end{document}